\documentclass[notitlepage,11pt]{article}
\usepackage{amssymb,amsmath, cases
}
\catcode`\@=11
\@addtoreset{equation}{section}

\catcode`\@=12
\usepackage{colortbl}%

\def\R{\mathbb{R}}

\def\f{\varphi}
\def\ep{\varepsilon}
\def\ep{\varepsilon}

\def\irn{\int\limits_{\R^n}}

\def\sstar{{2_m^*}}

\def\Dm{\left(-\Delta\right)^{\!m}\!}
\def\Dmhalf{\left(-\Delta\right)^{\!\frac{m}{2}}\!}
\def\Dshalf{\left(-\Delta\right)^{\!\frac{s}{2}}\!}
\def\Ds{\left(-\Delta\right)^{\!s}\!}

\def\tildeSobolev{\widetilde{H}^m(\Omega)}
\def\Hm{\widetilde{H}^m(\Omega)}							
\def\Hs{\widetilde{H}^s(\Omega)}							
\def\Dspace{\mathcal D^m(\R^n)}
\def\proof{\noindent{\textbf{Proof. }}}
\def\QED{\hfill {$\square$}\goodbreak \medskip}

\newtheorem{Theorem}{Theorem}[section]
\newtheorem{Lemma}[Theorem]{Lemma}
\newtheorem{Proposition}[Theorem]{Proposition}

\newtheorem{Remark}[Theorem]{Remark}
\newtheorem{Definition}[Theorem]{Definition}

\begin{document}

\title 
{Non-critical dimensions for critical problems\\
involving fractional Laplacians}

\author{Roberta Musina\footnote{Dipartimento di Matematica ed Informatica, Universit\`a di Udine,
via delle Scienze, 206 -- 33100 Udine, Italy. Email: {roberta.musina@uniud.it}. 
{Partially supported by Miur-PRIN 2009WRJ3W7-001 ``Fenomeni di concentrazione e {pro\-ble\-mi} di analisi geometrica''.}}~ and
Alexander I. Nazarov\footnote{
St.Petersburg Department of Steklov Institute, Fontanka 27, St.Petersburg, 191023, Russia, 
and St.Petersburg State University, 
Universitetskii pr. 28, St.Petersburg, 198504, Russia. E-mail: al.il.nazarov@gmail.com.
Supported by RFBR grant 11-01-00825 and by
St.Petersburg University grant 6.38.670.2013.}
}

\date{}

\maketitle

\begin{abstract}
We study the Brezis--Nirenberg effect in
two families of noncompact boundary value problems involving Dirichlet-Laplacian 
of arbitrary real order $m>0$. 
\footnotesize

\medskip

\medskip

\noindent
\end{abstract}

\normalsize

\bigskip

\section{Introduction}
\label{S:Introduction}

Let $m,s$ be two given real numbers, with $0\le s<m<\frac{n}{2}$. Let
$\Omega\subset\R^n$ be a bounded and smooth domain in $\R^n$
and put
$$
\sstar=\frac{2n}{n-2m}.
$$
We study equations
\begin{gather}
\label{eq:problem2}
(-\Delta)^mu=\lambda(-\Delta)^s u+|u|^{\sstar-2}u\quad\textrm{in $\Omega$,}\\
\nonumber
~\\
\label{eq:problem1}
(-\Delta)^mu=\lambda|x|^{-2s} u+|u|^{\sstar-2}u\quad\textrm{in $\Omega$,}
\end{gather}
under suitably defined
Dirichlet boundary conditions. In dealing with equation (\ref{eq:problem1}) we 
always assume that $\Omega$ contains the origin. For the definition
of fractional Dirichlet--Laplace operators $(-\Delta)^m, (-\Delta)^s$ and
for the variational approach to (\ref{eq:problem2}), (\ref{eq:problem1}) we refer to the next section.

The celebrated paper \cite{BN} by Brezis and Nirenberg
was the inspiration for a fruitful line of research about the effect of lower
order perturbations in noncompact variational problems. They
took as model the  case $n>2$, $m=1$, $s=0$, that is,
 \begin{equation}
\label{eq:BN_problem}
-\Delta u=\lambda u+|u|^{\frac{4}{n-2}}u\quad\textrm{in $\Omega$,}\qquad
u=0\quad
\textrm{on $\partial\Omega$.}
\end{equation}
Brezis and Nirenberg pointed out a remarkable phenomenon that appears
for positive values of the parameter $\lambda$: they proved 
existence of a nontrivial solution for any small $\lambda>0$ if $n\ge 4$; in contrast, 
in the lowest dimension
$n=3$ non-existence phenomena for sufficiently small $\lambda>0$ can be observed.
For this reason, the dimension $n=3$ has been  named {\em critical}\footnote{\footnotesize{
compare with \cite{PuSe}, \cite{GGS}.}}
for problem (\ref{eq:BN_problem}).

Clearly, as  larger $s$ is, as  stronger the effects of the lower order perturbations are
expected in equations (\ref{eq:problem2}), (\ref{eq:problem1}).
We are interested in the following question: 
{\em Given $m<\frac{n}{2}$, how large must be $s$ 
in order to have the existence of a ground state solution, for any
 arbitrarily small  $\lambda>0$ ?} In case of an affirmative answer,
 we  say that  $n$ is \underline{not} a critical dimension.

We present our main result,  that holds for any dimension $n\ge 1$
(see Section \ref{S:proof_main} for a more precise statement).

\bigskip
\noindent{\bf  THEOREM.}
{
\em If $s\ge 2m-\frac{n}{2}$ then $n$ is not a critical dimension for the Dirichlet boundary
value problems
associated to equations (\ref{eq:problem2})
and (\ref{eq:problem1}).}

\medskip

We point out some particular cases that are included in this result. 

\begin{description}
\item$\bullet$
If $m$ is an integer and $s=m-1$, then at most the lowest dimension
$n=2m+1$ is critical.

\item$\bullet$
For any $n>2m$ there always exist lower order perturbations
of the type $|x|^{-2s}u$ and of the type $\Ds u$ such that $n$ is not a critical dimension.

%
\item$\bullet$
If $ m<1/4$ then no dimension is critical, for any choice of $s\in[0,m)$.
\end{description}

After \cite{BN}, a large number of papers have been focussed
on studying the effect of linear perturbations in noncompact variational problems
of the type (\ref{eq:problem2}). Most of these papers deal with $s=0$,
when the problems (\ref{eq:problem2}) and (\ref{eq:problem1}) coincide. 
Moreover, as far as we know, all of them consider 
either polyharmonic case $2\le m\in \mathbb N$, see for instance
\cite{PuSe}, \cite{EFJ}, \cite{BG1}, \cite{Gr}, \cite{Gaz98}, or the case $m\in(0,1)$,
see \cite{SV1}, \cite{SV2}. We cite also \cite{CM2}, where equation
(\ref{eq:problem2}) is studied in case $m=2$, $s=1$. 
Thus, our Theorem \ref{T:main} covers all earlier existence results. 

Finally, we mention \cite{BaCPS} (see also \cite{Tan})
 where equation (\ref{eq:problem2}) for the so-called Navier-Laplacian is studied 
 in case $m\in(0,1)$,  $s=0$. For a comparison between the Dirichlet
 and Navier Laplacians we refer to \cite{FL}.

\medskip

The paper is organized as follows. After introducing
some notation and preliminary facts in
 Section \ref{S:Preliminaries}, we provide the main estimates in Section
 \ref{S:estimates}. In Section \ref{S:proof_main} we prove Theorem 1 and 
 point out an existence result for the case $s< 2m-\frac{n}{2}$.
 

 

\section{Preliminaries}
\label{S:Preliminaries}
The fractional Laplacian $\Dm u$ of a function $u\in {\cal C}^\infty_0(\R^n)$ is defined via the Fourier transform
$$
{\cal F}[u](\xi)=\frac{1}{(2\pi)^{n/2}}\int\limits_{\mathbb R^n} e^{-i~\!\!\xi\cdot x}u(x)~\!dx
$$
by the identity
\begin{equation}
\label{eq:classicalF}
{\cal F}\left[{\Dm u}\right](\xi)=|\xi|^{2m}{\cal F}[u](\xi).
\end{equation}
In particular, Parseval's formula gives
$$
\int\limits_{\mathbb R^n}\Dm u\cdot u~\!dx=
\int\limits_{\mathbb R^n}|\Dmhalf u|^2~\!dx=
\int\limits_{\mathbb R^n}|\xi|^{2m}|{\cal F}[u]|^2~\!d\xi~\!.
$$
We recall the well known Sobolev inequality
\begin{equation}
\label{eq:Sobolev}
\int\limits_{\R^n}|\Dmhalf u|^2~\!dx\ge {\cal S}_m\bigg(\,\int\limits_{\R^n}|u|^\sstar~\!dx\bigg)^{2/\sstar},
\end{equation}
that holds for any $u\in  {\cal C}^\infty_0(\R^n)$ and $m<\frac n2$, see for example  \cite[2.8.1/15]{Tr}.

Let $\Dspace$ be the Hilbert space obtained
by completing ${\cal C}^\infty_0(\R^n)$ with  respect to the Gagliardo norm 
\begin{equation}
\label{eq:Gagliardo}
\|u\|_m^2 =
\int\limits_{\R^n}|\Dmhalf u|^2~\!dx.
\end{equation}
Thanks
to (\ref{eq:Sobolev}), the space
$\Dspace$ is continuously embedded into $L^\sstar(\R^n)$. The 
{\em best Sobolev constant} ${\cal S}_m$ was explicitly computed in \cite{CoTa}. Moreover,
 it has been proved
in \cite{CoTa} that ${\cal S}_m$ is attained in $\Dspace$ by a unique family of functions, all of them
being obtained from
\begin{equation}
\label{eq:AT}
\phi(x)=(1+|x|^2)^{\frac{2m-n}{2}}
\end{equation}
by translations, dilations in $\R^n$ and multiplication by constants.

Dilations play a crucial role in the problems under consideration. Notice that
for any $\omega\in {\cal C}^\infty_0(\R^n)$, $R>0$ it turns out that
\begin{eqnarray}
\label{eq:dilation}
\int\limits_{\R^n}|\xi|^{2m}|{\cal F}[\omega](\xi)|^2~\!d\xi&=&
R^{n-2m}\int\limits_{\mathbb R^n}|\xi|^{2m}|{\cal F}[\omega(R\cdot)](\xi)|^2~\! d\xi\\
\int\limits_{\R^n}|\omega|^\sstar~\!dx&=&
R^n\int\limits_{\mathbb R^n}|\omega(R\cdot)|^\sstar~\! dx~\!.
\nonumber
\end{eqnarray}
Finally, we point out that the Hardy inequality
\begin{equation}
\label{eq:Hardy}
\int\limits_{\R^n}|\Dmhalf u|^2~\!dx\ge {\cal H}_m\int\limits_{\R^n}|x|^{-2m}|u|^2~\!dx
\end{equation}
holds for any function $u\in  {\Dspace}$. The 
{\em best Hardy constant} ${\cal H}_m$ was explicitly computed in \cite{He}.

The natural ambient space to study the Dirichlet boundary value problems for (\ref{eq:problem2}),  (\ref{eq:problem1}) is
$$\tildeSobolev=\{u\in \Dspace\,:\,{\rm supp}\, u\subset\overline{\Omega}\},
$$
endowed with the norm $\|u\|_m$. 
By Theorem 4.3.2/1 \cite{Tr}, for $m+\frac{1}{2}\notin\mathbb{N}$ this space coincides with $H^m_0(\Omega)$ (that is the closure of
${\cal C}^{\infty}_0(\Omega)$ in $H^m(\Omega)$), while for $m+\frac{1}{2}\in\mathbb{N}$ one has 
$\tildeSobolev\subsetneq H^m_0(\Omega)$. Moreover, ${\cal C}^{\infty}_0(\Omega)$ is dense in 
$\tildeSobolev$. 
Clearly, if $m$ is an integer then
$\tildeSobolev$ is the standard Sobolev space of 
functions $u\in H^m(\Omega)$ such that $D^\alpha u=0$
for every multiindex $\alpha\in\mathbb N^n$ with $0\le|\alpha|<m$.

We agree that $(-\Delta)^0u=u$, $\widetilde H^0(\Omega)=L^2(\Omega)$, since 
(\ref{eq:Gagliardo}) reduces to the standard $L^2$ norm in case $m=0$.\medskip

We define (weak) solutions of the Dirichlet problems for (\ref{eq:problem2}), (\ref{eq:problem1}) 
as suitably normalized critical points of the functionals
\begin{gather}
\label{eq:functional1}
{\cal R}^\Omega_{\lambda,m,s}[u]=
\frac
{\displaystyle \int\limits_\Omega|\Dmhalf u|^2~\!dx-\lambda\int\limits_\Omega|\Dshalf u|^2~\!dx}
{\bigg(\displaystyle\int\limits_{\Omega}|u|^{\sstar} dx\bigg)
^{2/\sstar\vphantom{\displaystyle 2^2}}}\\
\label{eq:functional2}
\widetilde {\cal R}^\Omega_{\lambda,m,s}[u]=
\frac
{\displaystyle \int\limits_\Omega|\Dmhalf u|^2~\!dx-\lambda\int\limits_\Omega|x|^{-2s}| u|^2~\!dx}
{\bigg(\displaystyle\int\limits_{\Omega}|u|^{\sstar} dx\bigg)
^{2/\sstar\vphantom{\displaystyle 2^2}}}~\!,
\end{gather}
respectively. It is easy to see that both functionals (\ref{eq:functional1}), 
(\ref{eq:functional2}) are well defined on $\Hm\setminus\{0\}$.\medskip

%
We conclude this preliminary section with some embedding results.

\begin{Proposition}
\label{P:Poincare}
Let $m,s$ be given, with $0\le s<m<n/2$. 
\begin{description}
\item$i)$ The space $\Hm$ is compactly embedded into $\Hs$.  In particular the infima
\begin{equation}
\label{eq:Poincare}
\Lambda_1(m,s):=\inf_{\scriptstyle u\in \Hm\atop\scriptstyle u\ne 0}
\frac{\|u\|_m^2}{\|u\|_s^2}~~,\qquad
\widetilde\Lambda_1(m,s):=\inf_{\scriptstyle u\in \Hm\atop\scriptstyle u\ne 0}
\frac{\|u\|_m^2}{\||x|^{-s}u\|_0^2}
\end{equation}
are positive and achieved.
\item$ii)$   
$\displaystyle{\inf_{\scriptstyle u\in \Hm\atop\scriptstyle u\ne 0}
\frac{\|u\|_m^2}
{\|u\|^2_{L^{2^*_m}}}={\cal S}_m.}$
\end{description}
\end{Proposition}

Statement $i)$ is well known for $\Lambda_1(m,s)$ 
and follows from
(\ref{eq:Hardy}) for $\widetilde\Lambda_1(m,s)$. To check $ii)$, use the inclusion
$\Hm\hookrightarrow \Dspace$ and a rescaling argument. Clearly, the Sobolev
constant ${\cal S}_m$ is never achieved on $\Hm$.

\section{Main estimates}
\label{S:estimates}
Let $\phi$ be the extremal of the Sobolev inequality (\ref{eq:Sobolev}) given by 
(\ref{eq:AT}). In particular, it holds that
\begin{equation}
\label{eq:M}
M:=\int\limits_{\R^n}|\Dmhalf \phi|^2~\!dx={\cal S}_m
\Big(\int\limits_{\R^n}|\phi|^{\sstar}~dx\Big)^{2/\sstar}.
\end{equation}
Fix $\delta>0$ and a cutoff function $\f\in {\cal C}^\infty_0(\Omega)$,
such that $\f\equiv 1$ on the  ball $\{|x|<\delta\}$ and $\f\equiv 0$ outside
$\{|x|<2\delta\}$. If $\delta$ is sufficiently small, the function
$$
u_\ep(x):=\ep^{2m-n}\f(x)\phi\left({\frac{x}{\ep}}\right)=\f(x)\left(\ep^2+|x|^2\right)^{\frac{2m-n}{2}}
$$
has compact support in $\Omega$. 
Next we  define
$$
\begin{array}{ll}
A^\ep_m:=\displaystyle\int\limits_{\Omega}|\Dmhalf u_\ep|^2 dx~&\quad
A^\ep_s:=\displaystyle\int\limits_{\Omega}|\Dshalf u_\ep|^2 dx~\\
&\\
\widetilde A^\ep_s:=\displaystyle\int\limits_{\Omega}|x|^{-2s}|u_\ep|^2 dx~&\quad
B^\ep:=\displaystyle\int\limits_{\Omega}|u_\ep|^{\sstar} dx
\end{array}
$$
and we denote by $c$ any universal positive constant.

\begin{Lemma}
\label{L:estimates}
It holds that
\begin{subnumcases}
{\label{eq:estimate}}
\label{eq:Hm_estimate}
A^\ep_m\le\ep^{2m-n}\left(M+c\ep^{n-2m}\right)&{}\\
\label{eq:Hs_estimate+}
A^\ep_s, \widetilde A^\ep_s\ge c\ep^{4m-n-2s}&{}\text{if $s>2m-\frac{n}{2}$}\\
\label{eq:Hs_estimate0}
A^\ep_s, \widetilde A^\ep_s\ge c~\!|\log\ep| &{}\text{if $s=2m-\frac{n}{2}$}\\    
\label{eq:Lp_estimate}
B^\ep\ge\ep^{-n}\left(({M}{\cal S}_m^{-1})^{\sstar/2}-c\ep^n\right)~\!.&{}
\end{subnumcases}

\end{Lemma}

\bigskip

\noindent
{\bf Proof of (\ref{eq:Hm_estimate}).}
First of all, from  (\ref{eq:dilation}) we get 
\begin{equation}
\label{eq:Aalpha}
A_m^\ep =\ep^{2m-n}\irn|\xi|^{2m}\left|{\cal F}\left[\f(\ep~\!\cdot)\phi\right]\right|^2~\!d\xi.
\end{equation}
Thus
$$
\Gamma_{m}^\ep:=\ep^{n-2m}A^\ep_m-M=
\irn|\xi|^{2{m}}\left|{\cal F}\left[\f(\ep~\!\cdot)\phi\right]\right|^2~\!d\xi-
\irn|\xi|^{2{m}}\left|{\cal F}[\phi]\right|^2~\!d\xi.
$$
We need to prove that
\begin{equation}
\label{eq:Gamma}
|\Gamma^\ep_m|\le c\ep^{n-2m}.
\end{equation}
If $m\in\mathbb N$, the proof of (\ref{eq:Gamma})
has been carried out in \cite{BN}, \cite{Gaz98}. Here we limit ourselves to the more 
difficult case, namely,  when $m$ is not an integer. We denote by $k:=\lfloor m\rfloor\ge 0$ the integer part
of $m$, so that $m-k>0$. Then
\begin{multline*}
\Gamma_m^\ep=
\irn|\xi|^{2k}{\cal F}[U_-]\cdot |\xi|^{2(m-k)}\overline{{\cal F}[U_+]}~\!d\xi\\
=2^{2(m-k)+\frac n2}\,\frac {\Gamma(m-k+\frac{n}2)}{\Gamma(-(m-k))}\cdot
\irn (-\Delta)^k U_-(x)\cdot
V.P.\irn\underbrace{\frac {U_+(x)-U_+(y)}{|x-y|^{n+2(m-k)}}}_{\Psi(x,y)}~\!dy~\!dx,
\end{multline*}
where $U_{\pm}=\f(\ep~\!\cdot~\!)\phi\pm\phi$ (the last equality follows from 
\cite[Ch. 2, Sec. 3]{GSh}).

We split the interior integral as follows:
$$V.P.\int\limits_{\R^n}\Psi dy=
\underbrace{V.P.\!\!\!\int\limits_{|y-x|\le\frac {|x|}2\atop ~}\Psi dy}_{I_1}+
\underbrace{\int\limits_{|y-x|\ge\frac {|x|}2\atop |y|\le|x|}\Psi dy}_{I_2}+
\underbrace{\int\limits_{|y-x|\ge\frac {|x|}2\atop |y|\ge|x|}\Psi dy}_{I_3}. 
$$
We claim that $|I_j|\le c|x|^{2k-n}$ for  $j=1,2,3$. Indeed, the Lagrange formula gives
\begin{multline*}
|I_1|
\le \max\limits_{|y-x|\le\frac {|x|}2}|D^2U_+(y)|\cdot
\int\limits_{|z|\le\frac {|x|}2}\frac {dz}{|z|^{n+2(m-k)-2}}\\
\le c|x|^{-(n-2m+2)}\cdot |x|^{2-2(m-k)}=c|x|^{2k-n}.
\end{multline*}
As concerns the last two integrals we estimate
$$|I_2|\le\int\limits_{|y-x|\ge\frac {|x|}2\atop |y|\le|x|}\frac {c|y|^{-(n-2m)}}{|x-y|^{n+2(m-k)}}\,dy
\le |x|^{-(n+2(m-k))}\cdot c|x|^{2m}=c|x|^{2k-n}
$$
and finally
\begin{eqnarray*}
|I_3|\le\int\limits_{|y-x|\ge\frac {|x|}2\atop |y|\ge|x|}\frac {c|x|^{-(n-2m)}}{|x-y|^{n+2(m-k)}}\,dy
&\le& c|x|^{-(n-2m)}\cdot\int\limits_{|z|\ge\frac {|x|}2}\frac {dz}{|z|^{n+2(m-k)}}\\
&\le& c|x|^{-(n-2m)}\cdot |x|^{-2(m-k)}=c|x|^{2k-n},
\end{eqnarray*}
and the claim follows.
Now, since
$$|(-\Delta)^k U_-(x)|\le \frac {c}{|x|^{n-2(m-k)}}~\!\chi_{\{|x|\ge\delta/\ep\}}+
\frac {c\ep^{2k}}{|x|^{n-2m}}~\!\chi_{\{\delta/\ep\le|x|\le2\delta/\ep\}},
$$
we obtain
$$|\Gamma_m^\ep|\le c\int\limits_{|x|\ge\delta/\ep}\frac {dx}{|x|^{2n-2m}}+c\int\limits_{\delta/\ep\le|x|\le2\delta/\ep}
\frac {\ep^{2k}\,dx}{|x|^{2n-2(m+k)}}\le c\ep^{n-2m},
$$
that completes the proof of (\ref{eq:Gamma}) and of (\ref{eq:Hm_estimate}).

\bigskip

\noindent{\bf Proof of (\ref{eq:Hs_estimate+}) and (\ref{eq:Hs_estimate0}).}
We use the Hardy inequality (\ref{eq:Hardy}) to get 
\begin{eqnarray*}
A^\ep_s&\ge& c\widetilde A^\ep_s
\ge
c\ep^{4m-2s-n}\int\limits_{\R^n}|x|^{-2s}|\f(\ep~\!\cdot)\phi|^2 dx\\
&\ge&c\ep^{4m-2s-n}\int\limits_{|x|<{\delta}/{\ep}}\frac{dx}{|x|^{2s}(1+|x|^2)^{n-2m}}~\!.
\end{eqnarray*}
The last integral converges as $\ep\to 0$ if $s>2m-\frac{n}{2}$, and  diverges 
with speed $|\log\ep|$ if $s=2m-\frac{n}{2}$.

\bigskip

\noindent{\bf Proof of (\ref{eq:Lp_estimate}).}
For $\ep$ small enough we estimate by below
\begin{eqnarray*}
\irn|u_\ep|^{\sstar}&=&\ep^{-n}\irn|\f(\ep~\!\cdot)\phi|^{\sstar}~\!dx   
=\ep^{-n}\Big(\irn|\phi|^{\sstar}~\!dx-\int\limits_{|x|>\delta/\ep}|\f(\ep~\!\cdot)\phi|^{\sstar}~\!dx\Big)\\
&\ge&\ep^{-n}\Big(({M}{\cal S}_m^{-1})^{\sstar/2}
-c\int\limits_{|x|>{\delta}/{\ep}}|x|^{-2n}~\!dx\Big)\\
&=& \ep^{-n}(({M}{\cal S}_m^{-1})^{\sstar/2}-c\ep^{n})  
\end{eqnarray*}
and the Lemma is completely proved.
\QED

\section{Two noncompact minimization problems}
\label{S:proof_main}

In this section we deal with the minimization problems 
\begin{equation*}
{\cal S}^\Omega_\lambda(m,s)=
\inf_{\scriptstyle u\in \Hm\atop\scriptstyle u\ne 0}{\cal R}^\Omega_{\lambda,m,s}[u];\qquad
\widetilde {\cal S}^\Omega_\lambda(m,s)=
\inf_{\scriptstyle u\in \Hm\atop\scriptstyle u\ne 0}\widetilde {\cal R}^\Omega_{\lambda,m,s}[u]~\!,
\end{equation*}
where the functionals ${\cal R}$ and $\widetilde {\cal R}$ are introduced in (\ref{eq:functional1})
and (\ref{eq:functional2}), respectively.


\begin{Lemma}
\label{L:standard1}
The following facts hold for any $\lambda\in\R$:
\begin{description}
\item$~~i)$ ${\cal S}^\Omega_\lambda(m,s) \le {\cal S}_m$;
\item$~ii)$  If $\lambda\le 0$ then ${\cal S}^\Omega_\lambda(m,s)= {\cal S}_m$ and it is not achieved;
 \item$iii)$ If $0<{\cal S}^\Omega_\lambda(m,s)< {\cal S}_m$, then ${\cal S}^\Omega_\lambda(m,s)$ is achieved.
 \end{description}
 The same statements hold for $\widetilde {\cal S}^\Omega_\lambda(m,s)$ instead of ${\cal S}^\Omega_\lambda(m,s)$.
\end{Lemma}

\proof
The proof is nowdays standard, and is essentially due to Brezis and Nirenberg \cite{BN}. 
We sketch it for the infimum ${\cal S}^\Omega_\lambda(m,s)$, for the convenience of the reader.

Fix $\ep>0$ and take $u\in {\cal C}^\infty_0(\R^n)\setminus\{0\}$ such that
\begin{equation}
\label{eq:Sob_eps}
({\cal S}_m+\ep)\bigg(\,\displaystyle\int\limits_{\R^n}|u|^{\sstar} dx\bigg)^{2/\sstar}
\ge \int\limits_{\R^n}|\Dmhalf u|^2~dx.
\end{equation}
Let $R>0$ be large enough, so that $u_R(\cdot):=u(R~\!\!\cdot)\in  {\cal C}^\infty_0(\Omega)$.
Using (\ref{eq:dilation}) we get
\begin{eqnarray*}
{\cal S}^\Omega_\lambda(m,s)&\le& 
\frac{\|u\|_m^2-
\lambda R^{2(s-m)}\|u\|_s^2}
{\|u\|^2_{L^{2^*_m}}}\le ({\cal S}_m+\ep)\left(1+cR^{2(s-m)}\right)~\!,
\end{eqnarray*}
where $c$ depends only on $u$ and $\lambda$. Letting $R\to\infty$
we get ${\cal S}^\Omega_\lambda(m,s)\le ({\cal S}_m+\ep)$ for any $\ep>0$, 
and  $i)$ is proved.

Next, if $\lambda\le0$ then clearly ${\cal S}^\Omega_\lambda(m,s)={\cal S}_m$.
If $\lambda=0$ then ${\cal S}_m$ is not achieved. The more it is not achieved for $\lambda<0$,
and $ii)$ holds.

%
%

Finally, to prove $iii)$ take a minimizing sequence $u_h$. It is convenient
to normalize $u_h$ with respect to the $L^{\sstar}$-norm, so that
$$
\displaystyle \int\limits_\Omega|\Dmhalf u_h|^2~\!dx-\lambda\int\limits_\Omega|\Dshalf u_h|^2~\!dx
={\cal S}^\Omega_\lambda(m,s)+o(1).
$$
We can assume that $u_h\to u$ weakly in $\Hm$ and strongly in $\Hs$  by Proposition~\ref{P:Poincare}.
Since
\begin{eqnarray*}
\lambda\int\limits_\Omega|\Dshalf u|^2~\!dx&=&\lambda\int\limits_\Omega|\Dshalf u_h|^2~\!dx+o(1)\\
&=& \int\limits_\Omega|\Dmhalf u_h|^2~\!dx-{\cal S}^\Omega_\lambda(m,s)+o(1)\\
&\ge& ({\cal S}_m-{\cal S}^\Omega_\lambda(m,s))+o(1),
\end{eqnarray*}
then $u\neq 0$. 
By the Brezis--Lieb lemma we get
$$
1=\|u_h\|^{2^*_m}_{L^{2^*_m}}=\|u_h-u\|^{2^*_m}_{L^{2^*_m}}+\|u\|^{2^*_m}_{L^{2^*_m}}+o(1).
$$
Thus
\begin{eqnarray*}
{\cal S}^\Omega_\lambda(m,s)&=&
{\|u_h\|_m^2-\lambda \|u_h\|^2_s}+o(1)\\
~&&\\
&=&
\Big(\|u_h-u\|_m^2+\|u\|_m^2\Big)-\lambda \Big(\|u_h-u\|^2_s+\|u\|_s^2\Big)+o(1)\\
~&&\\
&=&
\frac{\Big(\|u_h-u\|_m^2-\lambda \|u_h-u\|^2_s\Big)+
\Big(\|u\|^2_m- \lambda \|u\|_s^2\Big)}
{\Big(\|u_h-u\|^{2^*_m}_{L^{2^*_m}}+\|u\|^{2^*_m}_{L^{2^*_m}}\Big)^{2/2^*_m}}+o(1)\\
~&&\\
&\ge&
{\cal S}^\Omega_\lambda(m,s)\cdot\frac{\xi_h^{2}+1}{(\xi_h^{2^*_m}+1)^{2/2^*_m}}+o(1),
\end{eqnarray*}
where we have set
$$
\xi_h:=\frac{\|u_h-u\|_{L^{2^*_m}}}{\|u\|_{L^{2^*_m}}}.
$$

Since  $2^*_m>2$, this implies that $\xi_h\to 0$, that is,
$u_h\to u$ in $L^{2^*_m}$  and hence
$u$ achieves 
${\cal S}^\Omega_\lambda(m,s)$.
\QED
We are in position to prove our existence result, that  includes the 
theorem already stated in the introduction.

\begin{Theorem}
\label{T:main}
Assume $s\ge 2m-\frac{n}{2}$. 

\begin{description}
\item$~i)$~ If $0<\lambda<\Lambda_1(m,s)$ then  ${\cal S}^\Omega_\lambda(m,s)$ is achieved and
(\ref{eq:problem2}) has a nontrivial solution in $\Hm$.
\item$ii)$~ If $0<\lambda<\widetilde\Lambda_1(m,s)$ then  $\widetilde {\cal S}^\Omega_\lambda(m,s)$ is achieved and
(\ref{eq:problem1}) has a nontrivial solution in $\Hm$.
\end{description}
\end{Theorem}

\proof
Since $0<\lambda<\Lambda_1(m,s)$ then ${\cal S}^\Omega_\lambda(m,s)$ is positive, by Proposition
\ref{P:Poincare}. The main estimates in Lemma \ref{L:estimates} 
 readily imply
${\cal S}^\Omega_\lambda(m,s)<{\cal S}_m$. By Lemma \ref{L:standard1}, ${\cal S}^\Omega_\lambda(m,s)$ is achieved
by a nontrivial $u\in \Hm$, that solves (\ref{eq:problem2}) after  multiplication by a suitable constant.
Thus   $i)$ is proved. For $ii)$ argue in the same way.
\QED

In the case $s< 2m-\frac{n}{2}$ the situation is more complicated. We limit ourselves to point out
the next simple existence result.

\begin{Theorem}
\label{T:critical}
Assume $s< 2m-\frac{n}{2}$. 
\begin{description}
\item$~i)$ There exists $\lambda^*\in[0,\Lambda_1(m,s))$ such that
the infimum ${\cal S}^\Omega_\lambda(m,s)$ is attained for any $\lambda\in(\lambda^*,\Lambda_1(m,s))$, and
hence (\ref{eq:problem2}) has a nontrivial solution.
\item$ii)$ There exists $\widetilde\lambda^*\in[0,\widetilde\Lambda_1(m,s))$ such that
the infimum $\widetilde{\cal S}^\Omega_\lambda(m,s)$ is attained for any $\lambda\in(\widetilde\lambda^*,\widetilde\Lambda_1(m,s))$, and
hence (\ref{eq:problem1}) has a nontrivial solution.
\end{description}
\end{Theorem}

\proof
Use Proposition \ref{P:Poincare} to find $\f_1\in\Hm$, $\f_1\neq 0$, such that
$$
\int\limits_{\Omega}|\Dmhalf \f_1|^2~\!dx=\Lambda_1(m,s)
\int\limits_{\Omega}|\Dshalf \f_1|^2~\!dx~\!.
$$ Then test
${\cal S}^\Omega_\lambda(m,s)$ with $\f_1$  to get the strict inequality 
${\cal S}^\Omega_\lambda(m,s)<{\cal S}_m$.
The first conclusion follows by Proposition \ref{P:Poincare} and Lemma \ref{L:standard1}. For
(\ref{eq:problem1}) argue similarily.
\QED


\footnotesize
\label{References}

\end{document}